\documentclass[]{article}
\usepackage{amssymb, amsmath}
\author{Dioubina Anna\\
Tel Aviv University, Israel\\
 e-mail address: annadi@math.tau.ac.il}
\title{Instability of the virtual solvability and the property of
being virtually torsion-free for quasi-isometric groups}
\date{}

\begin{document}
\maketitle

\newtheorem{rem}{Remark}
\newtheorem{definition}{Definition}
\newtheorem{example}{Example}
\newtheorem{corollary}{Corollary}
\newtheorem{proposition}{Proposition}
\newtheorem{lemma}{Lemma}
\newtheorem{theorem}{Theorem}
\newcommand{\z}{\mathbb Z}

Many properties of finitely generated
 groups turn out to
be geometric. That is, they are preserved by quasi-isometries.
For example, the properties of being virtually free, hyperbolic, amenable,
finitely generated, virtually nilpotent or
virtually abelian are geometric \cite{Gr1}, \cite{Gr2}, \cite{GH}.
In this paper we construct examples of groups showing that virtual
solvability and the property of being virtually torsion-free
 are not preserved by
bi-Lipschitz maps (and hence by quasi-isometries). It provides the
negative answers
for questions formulated in \cite{GH}.

All groups under consideration are assumed to be finitely generated.

If $A$ is a (finitely generated) group, $d_A$ is a right invariant
 word metric on $A$
corresponding
to some symmetric finite system of generators, $l_A(x)=d_A(x,e)$.

If $f(x):A\to B$ then by $f(x)^{-1}$ we denote a function such that
$f(x)f(x)^{-1}=e$ for any $x\in A$.
\newline
\newline{\bf Definition.}
{\it Wreath product} of groups $A$ and $B$ is a semidirect product of
$A$ and $\oplus_A B$, where $A$ acts upon $\oplus_A B$ by shifts:
if $a\in A$, $f\in \oplus_A B$($\Leftrightarrow$ $f:A\to B$ with a finite
support),
 then $f^a(x)=f(xa^{-1}), x\in A$. That is for $(a_i,f_i)\in A\wr B$
$$
(a_1, f_1(x))(a_2,f_2(x))=(a_1a_2, f_1(xa_2^{-1})f_2(x) ).
$$
Denote the wreath product by $A\wr B$.
\newline
\begin{lemma}
Let $A, B$ be bi-Lipschitz equivalent groups. Then for any group $C$ the
groups
$C\wr A$ and $C\wr B$ are bi-Lipschitz equivalent.
\end{lemma}
{\bf Proof.}
Since $A$ and $B$ are bi-Lipschitz equivalent, there exists
$\varphi:A\to B$ such that for any $a\in A$
$$
K_2l_B(\varphi(a))\le l_A(a)\le K_1l_B(\varphi(a)).
$$
Define $\varphi_C :C\wr A\to C\wr B$ by
$$
\varphi_C( (c,f) )= (c,\varphi\circ f).
$$
Note that $\varphi_C$ is a one-to-one correspondence between $C\wr A$ and
$C\wr B$.
Let us show that $\varphi_C$ is bi-Lipschitz.
For any $\ae_1, \ae_2\in C\wr A: \ae_1=(c_1,f_1(x)), \ae_2=(c_2,f_2(x))$

$$
d_A(\ae_1,\ae_2)=l_A\left(\ae_1\ae_2^{-1}\right)=
l_A\left((c_1, f_1(x))(c_2,f_2(x))^{-1}\right)=
$$
$$
l_A\left(c_1c_2^{-1}, f_1(xc_2)f_2(xc_2)^{-1} \right).
$$
Similarly
$$
d_B(\varphi(\ae_1),\varphi(\ae_2))
=l_B(\varphi(\ae_1)\varphi(\ae_2)^{-1})=
l_B\left((c_1, \varphi\circ f_1(x))(c_2,(\varphi\circ f_2(x)))^{-1}\right)=
$$
$$
l_B\left(c_1c_2^{-1}, \varphi(f_1(xc_2))\varphi(f_2(xc_2)^{-1}) \right).
$$

Let
$$
 u(x)=f_1(xc_2)f_2(xc_2)^{-1},
$$
$$
w(x)=\varphi\left(f_1(xc_2))\varphi(f_2(xc_2)^{-1}\right).
$$
Note that for any $c\in C$
$$
K_2l_A(u(c))
 \le l_B(w(e))\le K_1l(u(e)),
$$
since
$$
l_A(u(c))=d_A(f_1(cc_2), f_2(cc_2))
$$
and
$$
l_B(w(c))=d_B(\varphi(f_1(cc_2)), \varphi(f_2(cc_2))).
$$
Let $a_1,a_2..., a_p$, $c_1,..., c_r$ be  sets of generators of $A$ and
$C$ respectively.
Let $\delta_x:C\to A$ be a function such that $\delta_x(e)=x$ and
$\delta_x(c)=e$ for $c\ne e$.
Consider
generators of $C\wr A$
$$
(c_1,\delta_e), ...,(c_r, \delta_e), (e,\delta_{a_1}), ..., (e,\delta_{a_p})
$$
and inverse ones.
Consider an analogous system of generators
 of $C\wr B$.

Note that for $c_0\in C$ and $f:C\to A$  with a finite support
$$
(c_i,\delta_e)(c_0,f)= (c_ic_0,f), (c_i,\delta_e)^{-1}(c_0,f)=(c_i^{-1}c_0,f)
$$
and
$$
(e,\delta_{a_i}(x))(c_0,f(x))=(c_0, \delta_{a_i}(xc_0^{-1})f(x)),
$$
$$
(e,\delta_{a_i}(x))^{-1}(c_0,f(x))=(c_0, \delta_{a_i^{-1}}(xc_0^{-1})f(x)).
$$
(In particular, multipliing by one of the generators either changes $c_0$
or the value of $f$ at $c_0$.)

Hence
$$
l_{C\wr A}(c_0,f)=K+\sum_{c\in C}l_A(f(c),
$$
where $K\ge 0$ depends on the support of $f$ and on $c_0$, but
 doesn't depend on $A$.

Then for any $\ae_1,\ae_2\in C\wr A$
$$
K_2 d_{C\wr A}(\ae_1,\ae_2)\le
d_{C\wr B}\left(\varphi_C(\ae_1), \varphi_C(\ae_2)\right)
\le K_1 d_{C\wr A}(\ae_1,\ae_2).
$$
Hence $\varphi_C$ is bi-Lipschitz.

\begin{proposition}
\begin{enumerate}
\item There exist quasi-isometric groups $G$ and $H$ such that $G$ is solvable,
but $H$ is not virtually solvable.
\item There exist quasi-isometric groups $G$ and $H$ such that $G$
has no torsion,
but no subgroup of finite index in $H$ is free of torsion.
\end{enumerate}
\end{proposition}
{\bf Proof.}
Let $D$ be a finite unsolvable group.
Consider $A=\z$, $B=Z\oplus D$, $C=\z$.
Set $G=C\wr A$ and $H=C\wr B$.
Then $G$ and $H$ are bi-Lipschitz equivalent, since  $A$ and $B$
are bi-Lipschitz equivalent.
$G$ is solvable (of degree 2) and $G$ has no torsion,
 but $H$ contains $\oplus_{\z} D$
 as a subgroup and hence is neither virtually solvable nor virtually
torsion-free.
In fact, all the elements of
$\oplus_{\z} D$ are torsion elements.
Any subgroup of finite index in $\oplus_{\z} D$ is not solvable,
since for some $i$ the projection
$$\pi_i:\oplus_{\z} D \to D$$
on the $i$-th component
is surjective.
Taking a subgroup preserves
virtual solvability,  since if $X, Y$ are subgroups
of $Z$, and $X$ is of finite index in $Z$, then $X\cap Y$ is of finite index
in $Y$.
Moreover, any subgroup of
finite index in $H$ has elements of order $|D|$.
\begin{rem}
There exist $G$ and $H$, $G$ is solvable, $H$ is not virtually solvable,
such that
Cayley graphs of $G$ and $H$ are
isometric for some choice of generators.
\end{rem}
{\bf Proof.}
Let $A, B$ be any finite groups such that $A$ is solvable,
$B$ is not solvable and $|A|=|B|$. Take all elements
as a system of generators of $A$ and $B$. Then Cayley graphs of $A$ and $B$
 are isometric.
Let $G=\z\wr A$, $H=\z\wr B$. Then Cayley graphs of $G$ and
$H$ coincide. The proof is analogous to the proof of Lemma 1.
It suffices
to note that Cayley graph of wreath product $X\wr Y$
depends only on $X$ and Cayley graph of $Y$ (in fact, it depends only
on Cayley graphs of $X$ and $Y$).
$A$ is solvable, but $B$ is not virtually solvable since it contains
$\oplus_{\z} B$ as a subgroup.

\end{document}